\documentclass[12pt]{article}
\usepackage{amssymb}
\def\overset#1#2{{\mathop{\kern0pt#2}\limits^{#1}}}
\begin{document}
\begin{center}
\bf PRODUCT PRESERVING BUNDLE FUNCTORS
ON MULTIFIBERED AND MULTIFOLIATE MANIFOLDS
\end{center}
%
%  Short title:
%
%  PRODUCT PRESERVING BUNDLE FUNCTORS
%
\begin{center}
\sc
Vadim V. Shurygin, jr.
\end{center}

\bigskip

\begin{center}
\parbox{0.8\textwidth}{\small
{\sc Abstract.} 
We show that the set of the equivalence classes of multifoliate structures is in one-to-one correspondence with the set of equivalence classes of finite complete projective systems of vector space epimorphisms.
After that we give the complete description of all product
preserving bundle functors on the categories of multifibered and
multifoliate manifolds.}
\end{center}

%\footnote{2000 {\it Mathematics
%Subject Classification}: 58A05, 58A32.\\
%{\it Key words and phrases}: product preserving bundle functors,
%multifoliate structures.}

\bigskip
\bigskip

2000 {\it Mathematics
Subject Classification}: 58A05, 58A32.

{\it Key words and phrases}: product preserving bundle functors,
multifoliate structures, projective systems.

%%
%%
%%     chapter 0
%%
%%

\bigskip
\bigskip

  In the middle 1980s Eck~\cite{Eck}, Kainz and Michor~\cite{Kainz-M}, 
Luciano~\cite{Luc} described all product preserving bundle functors on the
category of smooth manifolds in terms of  Weil bundles~\cite{Weil}
(see also~\cite{KMS}).
  In 1996 Mikulski~\cite{Mik} 
classified all product
preserving bundle functors on fibered manifolds.
  In the recent years Weil functors and product preserving functors
are of great interest, see e.g. Kol\'a\v r and Mikulski~\cite{K-M},
Kriegl and Michor~\cite{Kri-Mich}, Mu\~nos, Rodrigues, and Muriel~\cite{M-R-M}, Mikulski and Tom\'a\v s~\cite{MT, T}.
\par
  Kodaira and Spencer in~\cite{K-S} introduced the notion of a
multifoliate structure on  a smooth manifold.
  In the present paper, we introduce the category of multifibered
manifolds which is a subcategory of the category of multifoliate
manifolds and, following the lines of Mikulski~\cite{Mik},
describe all product preserving bundle functors on these
categories.
\par
  We denote the category of smooth manifolds by ${\cal M}f$
and that of fibered manifolds by ${\cal FM}$~\cite{KMS}.
  All manifolds and maps between manifolds under consideration are
assumed to be of class ${C}^\infty$.

\par\smallskip

%%
%%
%%     Section 1
%%
%%

\section{Projective systems of vector spaces}

  Let $(\Lambda=\{\alpha,\beta,\dots\}, \le)$
be a partially ordered set.
  A {\it projective system (an inverse system) over $\Lambda$}~\cite{A-M} is a collection
$(S_\alpha,\zeta_\alpha^\beta,\Lambda)$ consisting of sets
$S_\alpha$, $\alpha\in\Lambda$, and maps
$\zeta^\beta_\alpha:S_\beta\to S_\alpha$, $\alpha\le\beta$, called
projections, such that $\zeta^\alpha_\alpha={\rm id}_{S_\alpha}$
for all $\alpha\in\Lambda$ and
$\zeta^\beta_\alpha\circ\zeta_\beta^\gamma=\zeta_\alpha^\gamma$
when $\alpha\le\beta\le\gamma$.
  The projective limit of a projective system
$(S_\alpha,\zeta_\alpha^\beta,\Lambda)$ is the subset 
$$
S=
{\mathop{\null\rm lim}\limits_{\longleftarrow}}\,S_{\alpha}
\subset \prod_{\alpha\in\Lambda} S_\alpha
$$ 
consisting of all
elements $x=(x_\alpha)$ such that $\zeta^\beta_\alpha (x_\beta) =
x_\alpha$. If the set $S$ is not empty, then by $\zeta_\beta:S\to
S_\beta$ we denote the map which sends $x=(x_\alpha)$ into
$x_\beta$.
These maps are called {\it canonical projections}.
\par
  It will be convenient to denote projective
systems under consideration as follows:
$\zeta=(S_\alpha,\zeta_\alpha^\beta,\Lambda,S)$.
\par
  In this section, we will consider projective systems of vector spaces
$\xi=(L_\alpha,\xi_\alpha^\beta,\Lambda,L)$ satisfying the
following conditions:
\par
i) $L_\alpha$, $\alpha\in\Lambda$, and $L$
are  finite-dimensional vector spaces over $\mathbb R$;
\par
ii) all the maps $\xi^\beta_\alpha$ and $\xi_\alpha$ are linear
epimorphisms.
\par
\medskip
  By an isomorphism between two projective systems
$\xi=(L_\alpha,\xi_\alpha^\beta,\Lambda,L)$ and
$\xi'=(L_{\alpha'}',\xi_{\alpha'}^{\beta'},\Lambda',L')$ we  mean
a collection $(\omega, \{\psi_\alpha\}_{\alpha\in\Lambda})$
consisting of an isomorphism $\omega:\Lambda\to \Lambda'$ of
partially ordered sets and linear isomorphisms $\psi_\alpha :
L_\alpha \to L'_{\omega(\alpha)}$ such that
$\xi_{\omega(\alpha)}^{\omega(\beta)}\circ\psi_\beta =
\psi_\alpha\circ \xi_\alpha^\beta$ for $\alpha\le\beta$.
  An isomorphism $(\omega, \{\psi_\alpha\}_{\alpha\in\Lambda})$
gives rise to  the  isomorphism $\psi:=
{\mathop{\null\rm lim}\limits_{\longleftarrow}}\,\psi_{\alpha}:
L \to L'$  defined by
$\psi((x_\alpha)) = (\psi_\alpha(x_\alpha))$.
  The map $\psi$ is the unique isomorphism between $L$ and
$L'$ such that $\xi_{\omega(\alpha)} \circ \psi=\psi_\alpha \circ
\xi_\alpha$.
  Projective systems
$\xi=(L_\alpha,\xi_\alpha^\beta,\Lambda,L)$ and
$\xi'=(L_{\alpha'}',\xi_{\alpha'}^{\beta'},\Lambda',L')$ are said
to be {\it isomorphic} if there exists an isomorphism between
them.
\par
  An isomorphism from $\xi$ to itself
of the form $({\rm id}, \{f_\alpha\}_{\alpha\in\Lambda})$
is said to be {\it an automorphism of~$\xi$.}
  Denote by $GL(\xi)$ the group of all linear automorphisms
of~$L$ of the form
$f={\mathop{\null\rm lim}\limits_{\longleftarrow}}\,f_{\alpha}$
where $({\rm id}, \{f_\alpha\})$ is an automorphism of $\xi$.
\par\medskip
\noindent
{\bf Definition.}
A vector subspace  $K\subset L$
is said to be {\it invariant}
if every  $f\in GL(\xi)$ maps $K$ into itself.
\par\medskip
  One can easily see that the sum and the intersection of two invariant
subspaces are invariant subspaces. For any  $\alpha\in \Lambda$,
the subspace $K_\alpha:={\rm ker}\, \xi_\alpha \subset L$ is
invariant and $L_\alpha\cong L/K_\alpha$.
  In what follows we  will identify $L_\alpha$ and $L/K_\alpha$.
\par\medskip
\noindent
{\bf Definition.}
  A projective system
$\xi=(L_\alpha,\xi_\alpha^\beta,\Lambda,L)$ is said to be {\it
complete} if any   finite-codimensional invariant subspace of $L$ is of the form 
${\rm ker}\, \xi_\alpha$ for some $\alpha\in\Lambda$.
\par
\medskip
  Let  $\xi=(L_\alpha,\xi_\alpha^\beta,\Lambda,L)$ be a
projective system (not necessarily complete).
  Consider the set  $\{K_a\}_{a\in\widetilde\Lambda}$ of all finite-codimensional invariant subspaces $K_a$ of $L$.
  For any two invariant subspaces $K_a$, $K_b$ such that
$K_a \supset K_b$, denote by $\xi^b_a: L_b = L/K_b \to L_a =
L/K_a$ the canonical epimorphism.
  Let us endow $\widetilde\Lambda$
with the partial order defined as follows: $a \le b$ if and only if
$K_a \supseteq K_b$.
  One can easily see that the collection
$\widetilde\xi=(L_a,\xi_a^b,\widetilde\Lambda, \widetilde{L}=
{\mathop{\null\rm lim}\limits_{\longleftarrow}}\, L_a)$ is a
complete projective system.
  We  call it  {\it the completion of~$\xi$.}
  Obviously, $\widetilde\xi$ is complete.
  Since, for any  $\alpha\in\Lambda$, the subspace
$K_\alpha={\rm ker}\,\xi_\alpha$ is invariant, one can consider
$\Lambda$ as a subset of $\widetilde\Lambda$.
\par\medskip
\noindent
{\bf Definition.}
A projective system $\xi=(L_\alpha,\xi_\alpha^\beta,\Lambda,{L})$
is called {\it finite} if $\Lambda$ is finite.
\par\medskip
Obviously, when $\xi$ is finite, its limit $L$ is a finite-dimensional vector space.
\par\medskip
\noindent
{\bf Proposition 1.1.}
{\it If $\xi=(L_\alpha,\xi_\alpha^\beta,\Lambda,{L})$ is a finite complete projective system then 

$(1)$ $\Lambda$ contains the greatest element $\varepsilon$;

$(2)$ $L$ is isomorphic to  $L_\varepsilon$.}
\par\medskip
\noindent
{\bf Proof.}
Indeed, the zero subspace is invariant and of finite codimension.
\par\medskip
\noindent
{\bf Definition.}
Two projective systems
$\xi=(L_\alpha,\xi_\alpha^\beta,\Lambda,{L})$ and
$\xi'=(L_\alpha',\xi_{\alpha'}^{\beta'},\Lambda',{L}')$ are said
to be {\it equivalent} if there exists an isomorphism $\varphi: L
\to  L'$ such that $\varphi \circ f\circ \varphi^{-1} \in
GL(\xi')$ for any $f\in GL(\xi)$ and 
$$
\Phi: f\in GL(\xi) \mapsto
\varphi \circ f\circ \varphi^{-1} \in GL(\xi')
$$ 
is a group
isomorphism.
\par\medskip
  One can easily see that isomorphic projective systems are equivalent.
\par\medskip
\noindent
  {\bf Proposition 1.2.}
{\it Let $\xi=(L_\alpha,\xi_\alpha^\beta,\Lambda,{L})$ be a
projective system and
$\widetilde\xi=(L_a,\xi_a^b,\widetilde\Lambda,\widetilde L)$ its
completion.
  Then  $\xi$ and $\widetilde\xi$ are equivalent.}
\par\medskip
\noindent
  {\bf Proof.}
  In fact, the maps
$$
\varphi=(\varphi_\alpha):\widetilde L\to  L, \qquad
\varphi_\alpha: \widetilde L \ni
x = (x_a)_{a\in\widetilde\Lambda} \mapsto x_\alpha \in L_\alpha
$$
and
$$
\psi =(\psi_a):  L \to \widetilde L, \qquad
\psi_a: L\ni x \mapsto
x+K_a\in L_a= L/K_a
%\eqno(1.2)
$$
are mutually inverse isomorphisms
which induce an isomorphism of the groups
$GL(\xi)$ and $GL(\widetilde \xi)$.
$\Box$
\par\medskip
  The proof of the following proposition is immediate.
\par\medskip
\noindent
  {\bf Proposition 1.3.}
{\it If complete projective systems
$\xi=(L_\alpha,\xi_\alpha^\beta,\Lambda,{L})$ and
$\xi'=(L_\alpha',\xi_{\alpha'}^{\beta'},\Lambda',{L}')$ are
equivalent, then $\xi$ is isomorphic to $\xi'$.}
\par  \medskip\noindent
{\bf Definition.}
  Let $\xi=(L_\alpha, \xi_\alpha^\beta, \Lambda,  L)$
be a  projective system.
  A local diffeomorphism
$\varphi: U\subset L \to V\subset L$
between two open subsets of $L$ is called
{\it a $\xi$-diffeomorphism\/} if for any $x\in U$
there exist an open subset
$W(x)\subset U$ and a system of diffeomorphisms
$\{\varphi_\alpha : \xi_\alpha(W) \to
\xi_\alpha(\varphi(W))\}_{\alpha\in \Lambda}$ such that
$\varphi_\alpha \circ \xi_\alpha = \xi_\alpha \circ \varphi$
for any $\alpha\in \Lambda$.
\par
\medskip
  Denote the pseudogroup of all $\xi$-diffeomorphisms by $\Gamma(\xi)$.
  The tangent map $d\varphi_x$ of any $\xi$-diffeomorphism $\varphi:U\to V$ 
at every point $x\in U$ can be
viewed as an element of $GL(\xi)$.
\par
\medskip\noindent
{\bf Definition.}
  A {\it $\xi$-structure}
on an $n$-dimensional smooth manifold $(n=\dim L)$ is a maximal atlas
compatible with the pseudogroup $\Gamma(\xi)$.
  A smooth manifold endowed with a $\xi$-structure
is called a {\it $\xi$-manifold.}
\par
\medskip\noindent
{\bf Definition.}
  Let  $\xi=(L_\alpha, \xi_\alpha^\beta, \Lambda,  L)$ and
$\xi'=(L'_\alpha, {\xi'}_\alpha^\beta, \Lambda, L')$ be two projective
systems over the same partially ordered set $\Lambda$.
  A smooth map $g: U\subset L \to V\subset L'$ is called a
{\it $\Lambda$-smooth map} if for any $x\in U$ there exist
an open subset $W(x)\subset U$ and
a system of smooth maps
$\{g_\alpha : \xi_\alpha(W) \to
\xi'_\alpha(g(W))\}_{\alpha\in \Lambda}$ such that
$g_\alpha \circ \xi_\alpha = \xi'_\alpha \circ g$
for any $\alpha\in \Lambda$.
\par
\medskip\noindent
{\bf Definition.}
  Let  $M$ be a $\xi$-manifold and $M'$  a $\xi'$-manifold.
  A smooth map $f: M\to M'$ between a $\xi$-manifold $M$
and a $\xi'$-manifold $M'$ is called
a {\it $\Lambda$-smooth map} if it is $\Lambda$-smooth in terms of
the atlases defining $\xi$- and $\xi'$-structures on these manifolds.
\par
\medskip
  For  a fixed finite partially ordered set $\Lambda$,
all $\xi$-manifolds for all projective systems $\xi$ over $\Lambda$
together with $\Lambda$-smooth maps as morphisms
form a subcategory  ${\cal M}f_{\rm proj}(\Lambda)$
of the category   ${\cal M}f$.
\par\medskip
\noindent
  {\bf Proposition 1.4.}
{\it The category ${\cal M}f_{\rm proj}(\Lambda)$ admits products.}
\par\medskip
\noindent
  {\bf Proof.}
  Let $M$ be a $\xi$-manifold and $M'$  a $\xi'$-manifold,
where $\xi=(L_\alpha, \xi_\alpha^\beta, \Lambda,  L)$ and
$\xi'=(L'_\alpha, {\xi'}_\alpha^\beta, \Lambda, L')$.
  Then $M\times M'$ is a ($\xi\times \xi'$)-manifold, where 
$\xi\times \xi'$  denotes the projective system
$(L_\alpha\times L_\alpha', \xi_\alpha^\beta\times
{\xi'}_\alpha^\beta, \Lambda, L\times L')$.
$\Box$

%%
%%
%%            Section 2
%%
%%

\section{Multifoliate manifolds}

Multifoliate structures on smooth manifolds
were introduced by K.~Kodaira and D.C.~Spen\-cer~\cite{K-S}
as follows.
\par
\medskip
\noindent
{\bf Definition.}
  A pair $(\Lambda,p)$ consisting of
a finite partially  ordered set $\Lambda$ and
a surjective map
$$
p:\{1,\dots, n\}\ni i \, \mapsto \, p(i)\in \Lambda
$$
is called
{\it a multifoliate structure on the set $\{1,\dots, n\}$.}
\par
\medskip
  Denote by  $GL(\Lambda,p)$  the group of all linear isomorphisms
$$
f:{\mathbb R}^n\ni (x^i)\mapsto (f^i_j x^j)\in {\mathbb R}^n
$$
satisfying the condition
$$
f^i_j = 0 \quad\hbox{if}\quad p(i)\not\ge p(j),
$$
and by   $\Gamma(\Lambda,p)$ the pseudogroup of all
local diffeomorphisms
$g : U\subset {\mathbb R}^n \to V\subset {\mathbb R}^n$ such that
$ dg_x\in GL(\Lambda,p)$ for all $x\in U$.
\par \medskip
\noindent
{\bf Definition.}
A {\it $(\Lambda,p)$-multifoliate structure}
on an $n$-dimensional smooth manifold is a maximal atlas
compatible with the pseudogroup $\Gamma(\Lambda,p)$.
  We  call the local coordinates determined by a chart of this
atlas  {\it adapted  coordinates.}
  A smooth manifold endowed with a $(\Lambda,p)$-multifoliate
structure is called a {\it $(\Lambda,p)$-multifoliate manifold.}
\par
\medskip
\noindent
{\bf Definition.}
  Let $M$ be a  $(\Lambda,p)$-multifoliate manifold and  $N$ be a
$(\Lambda,p')$-multifoliate manifold.
  A {\it $\Lambda$-multifoliate map\/} $f: M\to N$
is a smooth map, satisfying the condition
$$
\frac{\partial f^a}{\partial x^i} = 0 \quad \hbox{ if } \quad
p'(a)\not\ge p(i)
$$
in adapted coordinates. Clearly, this definition does not depend on
the choice of a local coordinate system.
\par
\medskip
  For  a fixed finite partially ordered set $\Lambda$,
all $(\Lambda,p)$-multifoliate manifolds for all surjective maps
$p : \{1,\dots,n\} \to \Lambda$ and for all natural numbers
$n\ge {\rm card}\,\Lambda$
together with $\Lambda$-multifoliate maps as morphisms
form a subcategory  ${\cal M}f_\Lambda$
of the category   ${\cal M}f$.
We call it {\it the category of multifoliate
manifolds over $\Lambda$.}
\par\medskip
\noindent
  {\bf Proposition 2.1.}
{\it The category ${\cal M}f_\Lambda$ admits products.}
\par\medskip
\noindent
  {\bf Proof.}
If $M_a$ is a $(\Lambda, p_a)$-multifoliate manifold,
$p_a:\{1,\dots, n_a\}\to \Lambda$, $a=1, 2$,
then the product $M_1\times M_2$ is a $(\Lambda,p)$-multifoliate
manifold, where $p:\{1,\dots, n_1+n_2\}\to \Lambda$ is defined by
$$
p(i)=
\left\{
\begin{array}{ll}
p_1(i),  & i\le n_1;  \\[3pt]
p_2(i-n_1), \quad & i>n_1.
\end{array}
\right.
$$
\mbox{ } \hfill $\Box$
\par\medskip
\noindent
{\bf Definition.}
  We say that two multifoliate structures  $(\Lambda,p)$ and
$(\Lambda',p')$ on the same set $\{1,\dots,n\}$
are {\it equivalent} if there exists a linear automorphism
$\varphi: {\mathbb R}^n \to {\mathbb R}^n$ such that
$\varphi \circ f\circ \varphi^{-1} \in GL(\Lambda',p')$ for any
$f\in GL(\Lambda,p)$ and
$$
\Phi: f\in  GL(\Lambda,p) \mapsto 
\varphi \circ f\circ \varphi^{-1} \in GL(\Lambda',p')
$$
is a group isomorphism.
\par
\medskip
  Clearly, if $(\Lambda,p)$ is a multifoliate structure on
$\{1,  \dots, n\}$, then, for each permutation $\sigma$
on $\{1,  \dots, n\}$, the multifoliate structure
$(\Lambda,p\circ \sigma)$ is equivalent to $(\Lambda,p)$.
\par
  For a multifoliate structure $(\Lambda,p)$ on $\{1,\dots,n\}$,
define the sets $H_\alpha=\big\{i \,\big|\, p(i)\le\alpha\big\}$,
$\alpha\in\Lambda$, and let $k(\alpha)={\rm card}\, H_\alpha$.
  The vector spaces
$$
L_\alpha= \big\{(x^{i_1},\dots,x^{i_{k(\alpha)}})
\,\, \big| \,\, x^{i_s}\in{\mathbb R}, \, i_s\in H_\alpha, \, s=1,\dots,
k(\alpha)\big\}
$$ 
and the natural epimorphisms 
${\rm pr}^\beta_\alpha:L_\beta\to L_\alpha$,  $\alpha\le\beta$, form a
projective system whose limit can be naturally identified with
${\mathbb R}^n$.
  Denote this system and its completion, respectively,
by~$\xi(\Lambda,p)$ and~$\widetilde{\xi}(\Lambda,p)$.
\par\medskip
\noindent
{\bf Theorem 2.1.} \cite{ShVJr}
{\it
The correspondence $(\Lambda,p)\mapsto\widetilde{\xi}(\Lambda,p)$
induces a bijection between the
equivalence classes of multifoliate structures
$(\Lambda,p)$ and the equivalence classes of finite
complete projective systems of vector space epimorphisms.}
\par\medskip
\noindent
{\bf Proof.}
    We give here a sketch of the proof and refer for details
to~\cite{ShVJr}.
\par
  Show first that the correspondence
$(\Lambda,p)\mapsto\widetilde{\xi}(\Lambda,p)$
induces a map from the set of equivalence classes of multifoliate
structures to the set of equivalence classes of finite
complete projective systems of vector space epimorphisms.
  By Propositions 1.2 and 1.3,
it suffices to show that the groups $GL(\Lambda,p)$
and $GL(\xi(\Lambda,p))$ are isomorphic.
  In fact, there is a natural isomorphism
$GL(\Lambda,p) \to GL(\xi(\Lambda,p))$ which assigns to $g\in
GL(\Lambda,p)$ a collection of maps $\{g_\alpha: L_\alpha\to
L_\alpha\}$ defined as follows: $g_\alpha(y_\alpha)={\rm
pr}_\alpha(g(y))$ where $y\in {\mathbb R}^n$ is such that ${\rm
pr}_\alpha(y)=y_\alpha$. %%
\par\smallskip
  To prove that the correspondence indicated in the theorem is
one-to-one, we need to pass to the dual inductive system~\cite{A-M}.
\par
  Let $\xi=(L_\alpha,\xi_\alpha^\beta,\Lambda,L)$
be a finite complete projective system and let $\varepsilon\in\Lambda$ be the greatest element.
  The dual spaces~$L_\alpha^*$ together with the dual
maps~$(\xi_\alpha^\beta)^*$ form an inductive
system~$\xi^*=(L_\alpha^*,(\xi_\alpha^\beta)^*, \Lambda)$.
  The existence of the  greatest element  implies that the inductive limit of $\xi^*$ exists  and can be identified with the dual space $L^*$.
  Under this identification the dual maps $\xi^*_\alpha : L_\alpha^* \to L^*$ are
the canonical maps of $\xi^*$.
  Obviously, all the maps $(\xi^\beta_\alpha)^*$ and $\xi_\alpha^*$ are monomorhisms.
  We will call the inductive system 
$\xi^*=(L_\alpha^*,(\xi_\alpha^\beta)^*, \Lambda, L^*)$ the {\it dual of~$\xi$}.
\par
  For any $f\in GL(\xi)$ and for each $\alpha\in\Lambda$,
we have $\xi_\alpha^* \circ f_\alpha^* = f^*\circ \xi_\alpha^*$.
\par
  Denote by  $GL(\xi^*)$ the group of all linear automorphisms
$h:L^* \to L^*$ which are the limits of inductive systems of
linear automorphisms $h_\alpha:L_\alpha^*\to  L_\alpha^*$.
  Since the maps $\xi_\alpha^*$ are monomorphisms,
it will be convenient to consider each
$L_\alpha^*$ as a subspace of $L^*$.
  Then $h_\alpha = h|L_\alpha^*$ or, in other words, $h$
maps  $L_\alpha^*$ into itself.
  The correspondence $f\mapsto f^*$ is an isomorphism of the groups
$GL(\xi)$ and $GL(\xi^*)$.
\par
  The dual system
$\xi^*$ is {\it  complete} in the sense that it contains all
subspaces which are invariant with respect to each
$f^*\in GL(\xi^*)$.
\par
  By a {\it chain} in $\xi^*$ we  mean
a sequence of embeddings
$$
L_{\alpha_k}^*
\stackrel{(\xi_{\alpha_{k-1}}^{\alpha_k})^*}{\longleftarrow}
L_{\alpha_{k-1}}^*
\stackrel{(\xi_{\alpha_{k-2}}^{\alpha_{k-1}})^*}{\longleftarrow}
\dots \stackrel{(\xi^{\alpha_2}_{\alpha_1})^*}{\longleftarrow}
L_{\alpha_1}^*
$$
such that
$\alpha_1 < \alpha_2 < \dots < \alpha_{k}$
and $\alpha_i$ is the successor of $\alpha_{i-1}$ in $\Lambda$,
$i=2,\dots,k$, (that is, $\alpha_{i-1}\le \beta\le \alpha_i$ implies that either $\beta=\alpha_{i-1}$ or $\beta=\alpha_i$).
  The space  $L^*_{\alpha_{k}}$ is called the {\it end} of the chain.
\par
  $L_\alpha^*$ is said to be {\it a subspace of the first floor}
if $\alpha$ is a minimal element of~$\Lambda$.
  $L_\alpha^*$ is said to be  {\it a subspace of the
$k$-th floor} if each chain with end~$L^*_\alpha$ is of length no
greater than~$k$ and among all such chains there is at least one
of length $k$.
\par
  If $L_\alpha^*$ is a subspace of the first floor, we choose a basis
$B_\alpha=\{e_\alpha^1, \dots, e_\alpha^{s(\alpha)}\}$
in $L_\alpha^*$ and call the index
$\alpha$ {\it distinguished}.
  Let ${C}_1$ be the union of $B_\alpha$ for all subspaces of the
first floor.
  One can verify that the system ${C}_1$ is linearly independent.
  In fact, the assumption that the system is linearly dependent
contradicts the completeness of  $\xi^*$
(see \cite{ShVJr} for details).
\par
  Let now $L_\beta^*$ be a space of the second floor.
  Then  either  $L_\beta^* \subset {\cal L}\{{C}_1\}$,
where ${\cal L}\{{C}_1\}$  is the linear span of the system
${C}_1$, or one can choose a system of linearly independent elements
$B_\beta=\{e_\beta^1, \dots, e_\beta^{s(\beta)}\}$
in $L_\beta^*$ such that
$L_\beta^* = {\cal L}\{e_\beta^1, \dots, e_\beta^{s(\beta)}\}
\oplus (L_\beta^* \cap {\cal L}\{{C}_1\})$.
  In the latter case the index
$\beta$ is also  called {\it distinguished.}
  Let ${C}_2$ be the  union of $B_\beta$
for all subspaces of the second floor.
  The system  ${C}_1\cup {C}_2$
is linearly independent (see \cite{ShVJr} for details).
\par
  Suppose  that we have chosen systems
${C}_\ell$ for every $\ell\le k$.
  If $L_\gamma^*$ is a space of   $(k+1)$-th floor, then
either  $L_\gamma^*\subset{\cal L}_k :=
{\cal L} \{ {C}_1\cup\dots\cup {C}_k \}$ or
there exists a system of linearly independent  elements
$B_\gamma=\{e_\gamma^1$, \dots, $e_\gamma^{s(\gamma)}\}$
such that $L_\gamma^*=(L_\gamma^*\cap{\cal L}_k)\oplus
{\cal L}\{e_\gamma^1, \dots, e_\gamma^{s(\gamma)}\}$.
  In the latter case the index  $\gamma$
is called {\it distinguished}.
  Let ${C}_{k+1}$ be the  union of $B_\gamma$
for all  subspaces of the $(k+1)$-st floor.
  As above, the system ${C}_1\cup {C}_2\cup
\dots \cup{C}_{k+1}$ is linearly independent.
\par
  This process stops when we reach $L^*$.
  As a result, we obtain a subset $\Lambda'\subset\Lambda$ consisting of
distinguished elements and the corresponding basis $\{e^1,\dots,e^n\}$
in~$L^*$.
  Let $p:\{1,\dots,n\}\to \Lambda'$ be the map defined as follows:
$p(m)=\alpha$ where  $\alpha$ is the minimal
distinguished element such that $e^m \in L^*_\alpha$.
  The pair  $(\Lambda',p)$ is a multifoliate structure on $\{1,\dots, n\}$
and the group $GL(\xi)$ is isomorphic to $GL(\Lambda',p)$.
$\Box$
\par\medskip
\noindent
{\bf Corollary 2.1.}
{\it For any finite partially ordered set $\Lambda$, the categories
${\cal M}f_{\rm proj}(\Lambda)$ and ${\cal M}f_\Lambda$ are
isomorphic.}
\par\medskip
\noindent
{\bf Corollary 2.2.}
 {\it Let  $(\Lambda,p)$ be a multifoliate
structure on $\{1,\dots,n\}$ and $\xi(\Lambda,p)$  the
corresponding projective system.
  Let $\widetilde\xi(\Lambda,p)=(\widetilde L_a,\xi^b_a,\widetilde\Lambda,
\widetilde L)$ be the completion of $\xi(\Lambda,p)$ and
$(\Lambda', p')$ the multifoliate structure on $\{1,\dots,n\}$
determined by $\widetilde\xi(\Lambda,p)$.
Then
\par
$(1)$
the partially ordered sets $\Lambda$ and $\Lambda'$ are
canonically isomorphic;
\par
$(2)$
 the multifoliate structures  $(\Lambda,p)$ and $(\Lambda',p')$
are equivalent.}
\par \medskip
\noindent 
{\bf Proof.} (1) Every invariant subspace of
$\widetilde\xi(\Lambda,p)$  is of the form
$$
\widetilde L_{(\alpha_1, \dots,
\alpha_k)}= L/({\rm ker}\, \xi_{\alpha_1}\cap \dots \cap {\rm
ker}\, \xi_{\alpha_k}),
$$
 where $\alpha_1, \dots, \alpha_k\in
\Lambda$ are pairwise incomparable. Thus, $\widetilde\Lambda$ is
isomorphic to the set of all finite collections of pairwise
incomparable elements $(\alpha_1,\dots,\alpha_k)$ endowed with the
partial order defined as follows: $(\alpha_1,\dots,\alpha_k) \le
(\beta_1,\dots,\beta_\ell)$ if and only if ${\rm ker}\,
\xi_{\beta_1}\cap \dots \cap {\rm ker}\, \xi_{\beta_\ell}
\subseteq {\rm ker}\, \xi_{\alpha_1}\cap \dots \cap {\rm ker}\,
\xi_{\alpha_k}$.
\par
  The index $(\alpha_1, \dots, \alpha_k)\in\widetilde\Lambda$ is
distinguished if and only if $k=1$, and so the set of
all distinguished elements is naturally isomorphic to $\Lambda$.
\par
(2)
From Theorem 2.1 it follows that $GL(\Lambda,p)\cong GL(\xi(\Lambda,p))$
and $GL(\widetilde\xi(\Lambda,p))\cong GL(\Lambda',p')$. The rest
of the proof follows from Proposition 1.2.
$\Box$
\par\medskip
\noindent
{\bf Corollary 2.3.}
{\it
If $(\Lambda,p)$ and $(\Omega,q)$ are equivalent multifoliate
structures on $\{1,\dots,n\}$, then
\par
$(1)$ the partially ordered sets $\Lambda$  and
$\Omega$ are isomorphic;
\par
$(2)$
there exists a permutation $\sigma$ on $\{1,\dots,n\}$ such that
$q=p\circ \sigma$.}
\par \medskip
\noindent
{\bf Proof.}
(1)  Let $(\Lambda',p')$ and $(\Omega', q')$  be the
multifoliate structures corresponding to
the complete projective systems
$\widetilde\xi(\Lambda,p)$ and
$\widetilde\xi(\Omega,q)$ respectively.
  The systems $\widetilde\xi(\Lambda,p)$ and
$\widetilde\xi(\Omega,q)$ are equivalent.
  By Proposition 1.3, these systems are isomorphic.
  Hence the sets $\Lambda'$ and $\Omega'$ of their distinguished elements
are isomorphic.
  By Corollary 2.2, $\Lambda$ and $\Omega$ are isomorphic.
\par
(2)
  Let $\omega: \Lambda\ni\alpha\mapsto
\omega(\alpha)\in\Omega$ be an isomorphism.
  Recall that, for any distinguished index
$\beta\in\Lambda'\cong\Lambda$,
$s(\beta)$ denotes the number
of linearly independent elements in the system
$B_\beta$ defined in the proof of Theorem 2.1.
  One can easily see that the cardinality of the subset
$p^{-1}(\beta)\subset \{1,\dots,n\}$ coincides with $s(\beta)$.
  Since the projective systems
$\widetilde\xi(\Lambda,p)$ and $\widetilde\xi(\Omega,q)$ are
isomorphic, for every distinguished index $\beta\in\Lambda'\cong\Lambda$,
the numbers $s(\beta)$ and $s(\omega(\beta))$  coincide.
  This means that $p^{-1}(\alpha)$ and
$q^{-1}(\omega(\alpha))$ have the same cardinality for any
$\alpha\in\Lambda$.
  From this observation it follows that one can find a permutation
$\sigma$ on $\{1,\dots,n\}$ such that $q=p\circ \sigma$.
  In general, such a permutation is not unique.
$\Box$

%%
%%
%%   Section 3
%%
%%

\section{Multifibered manifolds. The
classification\protect\\ theorem}

\par\medskip
\noindent
{\bf Definition.}
  Let $\xi=(L_\alpha, \xi^\beta_\alpha, \Lambda, L)$ be a 
projective system of vector spaces
and let $\pi=(M_\alpha, \pi^\beta_\alpha, \Lambda, M)$  be a
projective system such that all $M_\alpha$ and
$M={\mathop{\null\rm lim} \limits_{\longleftarrow}}\,M_{\alpha}$
are smooth manifolds  and all
maps $\pi^\beta_\alpha: M_\beta \to M_\alpha$
and $\pi_\alpha: M \to M_\alpha$
are surjective submersions.
  Let $\cal A$ be a  $\xi$-structure on~$M$.
  We call 
$\pi=(M_\alpha, \pi^\beta_\alpha, \Lambda, M)$  a
{\it multifibered manifold\/} if the $\xi$-structure $\cal A$ on $M$ is compatible with all  projections $\pi_\alpha$ in the following sense:
for any point $x=(x_\alpha)\in M$, there are  charts $(U,h)$ centered at $x$ on $M$ and $(U_\alpha, h_\alpha)$ centered at $x_\alpha$ on $M_\alpha$, $\alpha\in \Lambda$, such that the following diagram commutes
$$
\begin{array}{rcl}
U\ &
\stackrel{\displaystyle h}{\hbox to 1.7cm{\rightarrowfill}}&
\ L \\[5pt]
\pi_\alpha \bigg\downarrow\
& &
 \,\,\bigg\downarrow \xi_\alpha \\
U_\alpha\, &
\stackrel{\displaystyle h_\alpha}{\hbox to 1.7cm{\rightarrowfill}}
& \ L_\alpha~.
\end{array}
$$

\par
\medskip

  It follows from  Corollary 2.1  that $M$
carries a structure of  multifoliate manifold.
  For any point $x=(x_\alpha)\in M$
the projective system of tangent spaces
$\xi_x=(T_{x_\alpha}M_\alpha, (d\pi^\beta_\alpha)_{x_\beta}, \Lambda, T_xM)$
is isomorphic to $\xi$.
\par
\medskip
\noindent
{\bf Definition.}
  {\it A multifibered map\/} $f:\pi\to\overline\pi$
between two multifibered manifolds
$\pi=(M_\alpha, \pi^\beta_\alpha, \Lambda,M)$ and
$\overline\pi=(\overline M_\alpha, \overline \pi^\beta_\alpha,\Lambda,
\overline{M})$
is a collection of maps
$\{f_\alpha: M_\alpha\to
\overline M_\alpha\}_{\alpha\in\Lambda}$ such that for all $\alpha\le\beta$
the diagram
$$
\begin{array}{rcl}
M_\beta\ &
\stackrel{\displaystyle f_\beta}{\hbox to 1.7cm{\rightarrowfill}}&
\ \overline M_\beta \\[5pt]
\pi^\beta_\alpha \bigg\downarrow\,\,\
& &
\ \,\,\bigg\downarrow \overline\pi^\beta_\alpha \\
M_\alpha\, &
\stackrel{\displaystyle f_\alpha}{\hbox to 1.7cm{\rightarrowfill}}
& \ \overline M_\alpha
\end{array}
$$
commutes.
  Each multifibered map  determines a unique smooth map
$f:M\to \overline M$.
\par\medskip
  Multifibered manifolds over $\Lambda$ together with multifibered maps
form a subcategory of the category ${\cal M}f_\Lambda$ of
multifoliate manifolds over $\Lambda$.
  We denote it by  ${\cal FM}_\Lambda$.
\par\medskip
\noindent
{\bf Proposition 3.1.}
{\it The category ${\cal FM}_\Lambda$ admits products.}
\par \medskip
\noindent
{\bf Proof.}
  If $\pi=(M_\alpha, \pi^\beta_\alpha, \Lambda,M)$ and
$\overline\pi=(\overline M_\alpha, \overline \pi^\beta_\alpha,
\Lambda,\overline{M})$
are two multifibered manifolds,
then their product is
$\pi\times\overline\pi:=(M_\alpha\times\overline M_\alpha,
\pi^\beta_\alpha\times\overline \pi^\beta_\alpha,\Lambda, M\times
\overline M)$.
$\Box$
\par\medskip
  The categories ${\cal FM}_\Lambda$ and ${\cal M}f_\Lambda$
are  local categories over manifolds.
\par
\medskip
\noindent
{\bf Definition.} An {\it inductive system of Weil algebra
homomorphisms over $\Lambda$} is a collection
$\mu=(A_\alpha,\mu_\alpha^\beta,\Lambda)$ consisting of Weil
algebras
$A_\alpha$, $\alpha\in\Lambda$, and Weil algebra homomorphisms
$\mu_\beta^\alpha:A_\alpha\to A_\beta$, $\alpha\le\beta$,
such that $\mu^\alpha_\alpha={\rm id}_{A_\alpha}$ for all
$\alpha\in\Lambda$ and
$\mu^\beta_\gamma\circ\mu_\beta^\alpha
=\mu^\alpha_\gamma$ when $\alpha\le\beta\le\gamma$.
  Let $\mu=(A_\alpha,\mu_\alpha^\beta,\Lambda)$ and
$\overline\mu=(\overline{A}_\alpha,
\overline{\mu}_\alpha^\beta,\Lambda)$ be two inductive systems of
Weil algebra homomorphisms.
  By a {\it morphism} $\nu :\mu\to \overline\mu$
we mean a  family
$\nu=(\nu_\alpha)_{\alpha\in\Lambda}$ of Weil algebra
homomorphisms
$\{\nu_\alpha: { A}_\alpha \to
\overline{ A}_\alpha\}$ such that for all $\alpha\le\beta$ the diagram
$$
\begin{array}{rcl}
{ A}_\alpha &
\stackrel{\displaystyle \nu_\alpha}{\hbox to 1.7cm{\rightarrowfill}}&
\ \overline{ A}_\alpha \\[5pt]
\mu^\alpha_\beta \bigg\downarrow\,\
& &
\ \,\,\bigg\downarrow \overline\mu^\alpha_\beta \\
{A}_\beta\, &
\stackrel{\displaystyle \nu_\beta}{\hbox to 1.7cm{\rightarrowfill}}
& \ \overline{A}_\beta
\end{array}
$$
commutes.

\par\medskip
\noindent
{\bf Theorem 3.1.} {\it Any  product preserving bundle functor
$F$ on the category
${\cal FM}_\Lambda$ or ${\cal M}f_\Lambda$
is uniquely determined by the inductive system
$\mu=(\mu^\alpha_\beta :{A}_\alpha \to
{A}_\beta)$ of Weil algebra homomorphisms. Any natural
transformation $\eta : F\to \overline F$ is uniquely determined by
the morphism $\nu :\mu\to \overline\mu$ of inductive systems of
Weil algebra homomorphisms.}
\par\medskip

  Since, by Theorem 2.1, any multifoliate manifold is locally
a multifibered manifold, it is enough to consider the case of a
bundle functor $F: {\cal FM}_\Lambda\to{\cal FM}$.
\par
   The proof of the Theorem 3.1 is essentially the same as the
Mikulski's proof
\cite{Mik} for the case of a bundle functor
${\cal FM}\to{\cal FM}$.
  We will reproduce the main scheme of the proof.
\par
\medskip
  Let   $\mu=(G_\alpha,\mu^\alpha_\beta,\Lambda)$
be an inductive system of
natural transformations of bundle functors, i.e.,
for any $\alpha\in\Lambda$, there is given a bundle functor
$G_\alpha: {\cal M}f\to {\cal FM}$ and for any $\alpha,\beta\in \Lambda$
such that $\alpha\le \beta$, there is given a natural transformation
$\mu^\alpha_\beta: G_\alpha\to G_\beta$ with the properties
$\mu^\alpha_\beta \circ \mu^\beta_\gamma =\mu^\alpha_\gamma$
and $\mu^\alpha_\alpha={\rm id}$.
  We define a bundle functor
$\prod\nolimits_\mu G_\alpha : {\cal FM}_\Lambda \to {\cal FM}$
as follows.
\par
  Consider a multifibered manifold
$\pi=(M_\alpha, \pi^\beta_\alpha, \Lambda,M)$.
  We let
$$
\prod\nolimits_\mu G_\alpha(\pi) \mathrel{:=}
\{(x_{\alpha}) \, | \,
G_\beta(\pi_\alpha^\beta) (x_{\beta}) =
\mu^\alpha_\beta(M_\alpha) (x_{\alpha})\}
\subset \prod\limits_{\alpha\in\Lambda} G_\alpha(M_\alpha).
$$
  The set $\prod\nolimits_\mu G_\alpha(\pi)$  is a submanifold
in $\prod_{\alpha\in\Lambda} G_\alpha(M_\alpha)$.
  We define the map
$$
p_\mu(\pi):\prod\nolimits_\mu G_\alpha(\pi) \to M
$$
as follows.
  Consider the bundle projection 
$$
\prod_{\alpha\in\Lambda} G_\alpha(M_\alpha) \, \to \,
\prod_{\alpha\in\Lambda} M_\alpha.
$$
  The image of its restriction to $\prod\nolimits_\mu G_\alpha(\pi)$ coincides
with
$M={\mathop{\null\rm lim}\limits_{\longleftarrow}}\,M_{\alpha}$,
thus defining the map
$p_\mu(\pi):\prod\nolimits_\mu G_\alpha(\pi) \to M$ which is a
surjective submersion.
\par\smallskip
  Let $f=(f_\alpha):\pi \to \overline\pi$
be a multifibered map.
  We set
$$
\prod\nolimits_\mu G_\alpha(f) \mathrel{:=}
\hbox{ the restriction of }
\prod_{\alpha\in\Lambda} G_\alpha(f_\alpha).
$$
  The map 
$$
\prod\nolimits_\mu G_\alpha(f) : \prod\nolimits_\mu G_\alpha(\pi) \to
\prod\nolimits_\mu G_\alpha(\overline{\pi}) 
$$
is well-defined
since all  $\mu^\alpha_\beta$ are natural transformations.
\par
  The correspondence
$$
\prod\nolimits_\mu G_\alpha:{\cal FM}_\Lambda\to{\cal FM}
$$
is a bundle functor.
\par
  Now let
$\overline\mu=(\overline \mu^\alpha_\beta:\overline G_\alpha
\to \overline G_\beta)$ be another inductive system of natural
transformations of bundle functors.
  Suppose that there is given a family
$\nu=(\nu_\alpha: G_\alpha\to\overline
G_\alpha)$ of natural transformations such that, for any manifold $X$
and $\alpha\le \beta$, the diagram
%$$
\begin{equation}
\begin{array}{rcl}
G_\alpha(X) &
\stackrel{\displaystyle \nu_\alpha(X)}{\hbox to 1.7cm{\rightarrowfill}}&
\ \overline G_\alpha(X) \\[5pt]
\mu_\beta^\alpha(X) \bigg\downarrow~~~
& &
~~~~\bigg\downarrow \overline\mu_\beta^\alpha(X) \\
G_\beta(X) &
\stackrel{\displaystyle \nu_\beta(X)}{\hbox to 1.7cm{\rightarrowfill}}
& \ \overline G_\beta(X)
\end{array}
\label{G-diag}
\end{equation}
%$$
commutes.
  Then we define the natural transformation
$$
\prod_{\mu,\overline\mu} \nu_\alpha:
\prod\nolimits_\mu G_\alpha \to
\prod\nolimits_{\overline\mu} \overline G_\alpha
$$
as follows.
\par
   For a multifibered manifold
$\pi=(M_\alpha, \pi^\beta_\alpha, \Lambda,M)$,
we define the map
$$
\prod\nolimits_{\mu,\overline\mu} \nu_\alpha (\pi) :
\prod\nolimits_\mu G_\alpha(\pi) \,\to\,
\prod\nolimits_{\overline\mu} \overline G_\alpha(\pi)
$$ 
to be  the restriction of 
$\prod_{\alpha\in\Lambda} \nu_\alpha(M_\alpha)$.
  Since each $\nu_\alpha$ is a  natural transformation, the map
$\prod_{\mu,\overline\mu} \nu_\alpha(\pi)$  is well-defined.
  The family  
$$
  \prod_{\mu,\overline\mu} \nu_\alpha=
\Bigl\{ \prod_{\mu,\overline\mu} \nu_\alpha(\pi)\Bigr\}:
\prod\nolimits_\mu G_\alpha \to
\prod\nolimits_{\overline\mu} \overline G_\alpha
$$
 is a natural
transformation.

\par\smallskip
  Let us denote by $pt$ a one-point manifold.
  Consider a smooth manifold $X$.
  For any  $\alpha\in\Lambda$, we construct a multifibered manifold
$i_\alpha(X)=(X_\gamma,r^\gamma_\delta,\Lambda,X)$ in the following way.
  We let  $X_\gamma=X$ if $\gamma\ge \alpha$, and $X_\gamma=pt$ otherwise.
  Each projection $r^\gamma_\delta$ is either the identity map
${\rm id}_X: X\to X$ if $\gamma\ge \alpha$, $\delta\ge \alpha$, or the unique map 
$pt_X: X\to pt$ if $\gamma\ge\alpha$, $\delta\not\ge\alpha$, or the unique map
$pt\to pt$ if $\gamma\not\ge\alpha$, $\delta\not\ge\alpha$.
  Clearly, ${\mathop{\null\rm lim}
\limits_{\longleftarrow}}\,X_{\alpha}=X$.
  We can consider any map $f: X\to Y$ as a multifibered map
$f: i_\alpha(X)\to i_\beta(Y)$ for $\alpha\le\beta$.
  Thus we obtain the bundle functors
$$
i_\alpha: {\cal M}f\to {\cal FM}_\Lambda
$$
and the natural
transformations
$$
{\rm id}^\alpha_\beta : i_\alpha\to i_\beta, \quad \alpha\le\beta ,
$$
consisting of 
${\cal FM}_\Lambda$-morphisms ${\rm id}_X: i_\alpha(X)\to i_\beta(X)$.
  Obviously, the functors $i_\alpha$ preserve products.
\par
  Let $F: {\cal FM}_\Lambda \to {\cal FM}$ be a bundle functor.
  Consider the bundle functors
%$$
\begin{equation}
G^F_\alpha
\mathrel{\colon=} F \circ i_\alpha : {\cal M}f \to
{\cal FM}.
\label{GFa}
\end{equation}
%$$
  If  $F$ preserves products, then  the functors $G^F_\alpha$ also preserve products.
\par
  We define an inductive system
$\mu^F=((\mu^F)^\alpha_\beta)$
of natural transformations  as follows:
%$$
\begin{equation}
(\mu^F)^\alpha_\beta \mathrel{:=} F({\rm id}^\alpha_\beta):
G^F_\alpha\to G^F_\beta.
\label{muF}
\end{equation}
%$$
\par
  Let $\overline F: {\cal FM}_\Lambda \to {\cal FM}$
be another bundle functor, and let
$\eta=\{\eta_\pi\}: F\to \overline F$ be a natural transformation.
  We define the family of natural transformations
$$
\nu^\eta = (\nu^\eta_\alpha : G^F_\alpha \to G^{\overline F}_\alpha)
$$
by
%$$
\begin{equation}
  \nu_\alpha^\eta (X) \mathrel{:=} \eta_{i_\alpha(X)} :
  G^F_\alpha(X) \to G^{\overline F}_\alpha(X)
\label{nu_eta}
\end{equation}
%$$
for any manifold  $X$.
  The diagram
$$
\begin{array}{rcl}
G^F_\alpha(X) &
\stackrel{\displaystyle \nu_\alpha^\eta(X)}{\hbox to 2cm{\rightarrowfill}}&
\ G^{\overline F}_\alpha(X) \\[5pt]
(\mu^F)_\beta^\alpha (X)\bigg\downarrow\ \quad
& &
\ \quad\bigg\downarrow (\mu^{\overline F})_\beta^\alpha(X) \\
G^F_\beta(X) &
\stackrel{\displaystyle \nu_\beta^\eta(X)}{\hbox to 2cm{\rightarrowfill}}
& \ G_\beta^{\overline F}(X)
\end{array}
$$
commutes  for any manifold $X$ and any $\alpha\le \beta$.
\par
\medskip
  Let  $F: {\cal FM}_\Lambda \to {\cal FM}$ be a bundle functor.
  Following Mikulski, we construct a natural transformation
$$
\Theta =\{\Theta_\pi\}: F \to
\prod\nolimits_{\mu^F} G_\alpha^F.
$$
\par
  Let  $\pi=(M_\alpha, \pi^\beta_\alpha, \Lambda,M)$
be a multifibered manifold.
  For any $\alpha\in\Lambda$ we define a multifibered map
$j_\alpha: \pi\to i_\alpha(M_\alpha)$
as follows: we let
$(j_\alpha)_\gamma := \pi^\gamma_\alpha$ if $\alpha\le\gamma$ and
$(j_\alpha)_\gamma := pt_{M_\gamma}$ otherwise.
\par
  The image of the map
$$
\prod_{\alpha\in\Lambda} F(j_\alpha) :F(\pi) \to
\prod_{\alpha\in\Lambda} F(i_\alpha(M_\alpha)) =
\prod_{\alpha\in\Lambda} G^F_\alpha(M_\alpha)
$$
is contained in $\prod\nolimits_{\mu^F} G_\alpha^F\,(\pi)$.
%
%  This follows  immediately  from the equality
%${\rm id}^\alpha_\beta(M_\alpha) \circ j_\alpha =
%i_\beta(\pi^\beta_\alpha) \circ j_\beta$.
%
   Therefore, the map
$$
  \Theta_\pi \mathrel{:=} \prod_{\alpha\in\Lambda}F(j_\alpha) :
F(\pi) \to
\prod\nolimits_{\mu^F} G_\alpha^F\,(\pi) \subset
\prod_{\alpha\in\Lambda} G^F_\alpha(M_\alpha).
$$
is well-defined.
 \par
  The family $\Theta = \{\Theta_\pi\} : F \to
\prod\nolimits_{\mu^F} G_\alpha^F$ is a natural transformation.
\par

%%%%%%%%%%%%%%%%%%%%%%%%%%%%%%%%%%%%%%%%%%%%%%%%%%

  Let now $\mu=(\mu^\alpha_\beta : G_\alpha\to G_\beta)$
be an inductive system of natural transformations of bundle functors
$G_\alpha: {\cal M}f\to {\cal FM}$.
   Consider the corresponding bundle functor
$F:=\prod\nolimits_\mu G_\alpha : {\cal FM}_\Lambda
\to{\cal FM}$.
  Denote by
$\overset\circ\mu=(\overset\circ\mu{}^\alpha_\beta:
\overset\circ{G}_\alpha\to \overset\circ{G}_\beta)$
the corresponding inductive system of natural transformations~(\ref{muF}).
  Then
$$
\overset\circ{G}_\gamma(X) =
\{ x_{\alpha} \in G_\alpha(X_\alpha) \;| \;
(\mu^F)^\alpha_\beta(M_\alpha) (x_{\alpha}) =
G^F_\beta(\pi_\alpha^\beta) (x_{\beta})\}
\subset \prod_{\alpha\in\Lambda} G_\alpha(X_\alpha),
$$
where $X_\alpha=X$ for $\alpha\ge \gamma$, otherwise $X_\alpha=pt$.
\par
  For any manifold $X$ and for any $\alpha\in\Lambda$, we define the map
%$$
\begin{equation}
{\cal O}_\alpha(X) : \overset\circ{G}_\alpha(X) \to G_\alpha(X)
\label{O_alpha}
\end{equation}
%$$
as the restriction of the standard projection
$\prod_{\beta\in\Lambda} G_\beta(X_\beta) \to G_\alpha(X)$.
\par
  The families 
$$
{\cal O}_\alpha=\{{\cal O}_\alpha(X) \} :
\overset\circ{G}_\alpha\to G_\alpha
$$
are natural transformations.
  They all are natural equivalences if and only if
every map $\mu^\alpha_\beta(pt):G_\alpha(pt)\to G_\beta(pt)$
is a diffeomorphism.
  The diagram
$$
\begin{array}{rcl}
\overset\circ{G}_\alpha(X) &
\stackrel{\displaystyle {\cal O}_\alpha(X)}{\hbox to 2cm{\rightarrowfill}}&
\ G_\alpha(X) \\[5pt]
\overset\circ\mu{}_\beta^\alpha (X)\bigg\downarrow\ \quad
& &
\ \quad\bigg\downarrow \mu_\beta^\alpha(X) \\
\overset\circ{G}_\beta(X) &
\stackrel{\displaystyle {\cal O}_\beta(X)}{\hbox to 2cm{\rightarrowfill}}
& \ G_\beta(X)
\end{array}
$$
is commutative for any manifold $X$ and any $\alpha\le\beta$.
\par
  Suppose now that
the inductive system $\mu=(\mu^\alpha_\beta : G_\alpha\to G_\beta)$
satisfies the condition that
all the maps
$\mu^\alpha_\beta(pt) : G_\alpha(pt)\to G_\beta(pt)$ are
diffeomorphisms.
  \par
  For any  multifibered manifold
$\pi=(M_\alpha, \pi^\beta_\alpha, \Lambda,M)$
we let
$$
T^\mu(\pi)=
\left\{
\begin{array}{ll}
G_\alpha(X) &
\hbox{ if } \pi=i_\alpha(X) \hbox{ for some }
X\in{\cal M}f,\ \alpha\in\Lambda,  \\
\prod\nolimits_\mu G_\alpha\, (\pi) &
\hbox{ otherwise.}
\end{array}
\right.
$$
  Then $T^\mu(\pi)$ is a fibered  manifold over
$M={\mathop{\null\rm lim}\limits_{\longleftarrow}}\,M_{\alpha}$.
  We also define the map 
$$
I_\pi : T^\mu(\pi) \to
\prod\nolimits_\mu G_\alpha\, (\pi)
$$
as follows:
$$
I_\pi=
\left\{
\begin{array}{ll}
{\cal O}_\alpha^{-1}(X) & \hbox{ if } \pi=i_\alpha(X), \\
{\rm id}_ {\Pi_\mu\> G_\alpha\, (\pi)}~ &
\hbox{ otherwise,}
\end{array}
\right.
$$
where  ${\cal O}_\alpha(X) : \overset\circ{G}_\alpha (X)=
\prod\nolimits_\mu G_\alpha\, (i_\alpha(X)) \to G_\alpha(X)$
are defined by~(\ref{O_alpha}).
  We let
$$
  T^\mu(f):=
  I_{\overline\pi}^{-1} \circ
    \prod\limits_{\alpha\in\Lambda} G_\alpha(f_\alpha)
       \circ I_\pi. 
$$
  The correspondence  $T^\mu$ thus defined is a bundle functor
${\cal FM}_\Lambda \to{\cal FM}$, and the family
$$
I=\{I_\pi\} : T^\mu\to \prod\nolimits_\mu G_\alpha
$$
is a natural transformation.
\par
\medskip
  If all $G_\alpha$  preserve products, then
$G_\alpha(pt)=pt$, hence the maps
$\mu^\alpha_\beta(pt)$ are diffeomorphisms.
  In this case, $I$ is a natural equivalence and the functor
 $T^\mu$ also preserves products.
\par
  Let now
$\overline\mu=(\overline \mu^\alpha_\beta: \overline G_\alpha
\to \overline G_\beta)$ be another inductive system of natural
transformations such that all
$\overline\mu^\alpha_\beta(pt)$ are diffeomorphisms,
and let
$\nu=(\nu_\alpha: G_\alpha\to\overline G_\alpha)$ 
be a family of natural transformations such that the diagram
(\ref{G-diag}) is commutative.
  Following Mikulski, we define a natural transformation
$\widetilde\nu=\{\widetilde\nu_\pi\} : T^\mu \to
T^{\overline\mu}$ to be the composition
$$
\widetilde\nu_\pi :
T^\mu(\pi)
\stackrel{\scriptstyle I_\pi}{\hbox to 1cm{\rightarrowfill}}
\prod\nolimits_{\mu} G_\alpha\, (\pi)
\stackrel{\scriptstyle \Pi_{\mu,\overline\mu}\>
\nu_\alpha(\pi)}{\hbox to 1.9cm{\rightarrowfill}}
\prod\nolimits_{\overline\mu} \overline G_\alpha\, (\pi)
\stackrel{\scriptstyle \overline I_\pi^{-1}}{\hbox to 1cm{\rightarrowfill}}
T^{\overline\mu}(\pi)
$$
for any multifibered manifold $\pi$.
\par
  In the case  $F=T^\mu$  the natural transformations
$(\mu^F)^\alpha_\beta : G^F_\alpha \to G^F_\beta$ coincide with
$\mu^\alpha_\beta$, i.e.,
$$
  \mu^F = \mu \quad \hbox{if}\quad F=T^\mu.
$$
  
\par
  Let $F:{\cal FM}_\Lambda\to {\cal FM}$ be a bundle functor such that
$(\mu^F)^\alpha_\beta(pt)$, $\alpha\le\beta$, are diffeomorphisms.

  Then we define a natural transformation
$\kappa=\{\kappa_\pi\} : F \to T^{\mu^F}$ to be the composition
%$$
\begin{equation}
 \kappa_\pi  : F(\pi)
\stackrel{\scriptstyle \Theta_\pi}{\hbox to 1cm{\rightarrowfill}}
\prod\nolimits_{\mu^F} G^F_\alpha\, (\pi)
\stackrel{\scriptstyle I_\pi^{-1}}{\hbox to 1cm{\rightarrowfill}}
T^{\mu^F}(\pi)
\label{kappa}
\end{equation}
%$$
for any multifibered manifold $\pi$.

The proofs of the following propositions are similar to the proofs of Theorems~2.1 and~2.2 in~\cite{Mik}.
\par
\medskip
\noindent
{\bf Proposition 3.2.}
{\it
 $(1)$ Let $F: {\cal FM}_\Lambda \to {\cal FM}$
be a product preserving  bundle functor.
Then the natural transformation  $\kappa : F\to
T^{\mu^F}$ is a natural equivalence.}
\par
{\it
  $(2)$ If $\mu=(\mu^\alpha_\beta:G_\alpha\to G_\beta)$
is an inductive system of natural transformations between product
preserving bundle fuctors $G_\alpha :{\cal M}f\to {\cal FM}$
and $\kappa$ is the natural transformation
$(\ref{kappa})$ for $F=T^\mu$, then $\kappa:T^\mu\to T^\mu$ and $\kappa_\pi =
{\rm id}_{T^\mu(\pi)}$ for any multifibered manifold $\pi$.}
\par
{\it
  $(3)$ For  $\mu=(\mu^\alpha_\beta:G_\alpha\to G_\beta)$
the functor $T^\mu$ is a
product preserving bundle functor on the category ${\cal FM}_\Lambda$
unique up to a natural equivalence
such that the natural transformation
$\mu^F$ corresponding to
$F=T^\mu$ coincides with~$\mu$.}
\par \medskip

\noindent
{\bf Proposition 3.3.}
{\it Let $F, \overline F : {\cal FM}_\Lambda\to{\cal FM}$
be two product preserving bundle functors.
Let $\mu^F=((\mu^F)^\alpha_\beta : G^F_\alpha\to G^F_\beta)$ and
$\mu^{\overline F}=((\mu^{\overline F})^\alpha_\beta:
G^{\overline F}_\alpha\to G^{\overline F}_\beta)$
be the corresponding inductive systems of natural transformations.
Let  $\nu=(\nu_\alpha: G^F_\alpha\to
G^{\overline F}_\alpha)$ be the family of natural transformations
such that the diagram
$$
\begin{array}{rcl}
G^F_\alpha(X) &
\stackrel{\displaystyle \nu_\alpha(X)}{\hbox to 2cm{\rightarrowfill}}&
\ G^{\overline F}_\alpha(X) \\[5pt]
(\mu^F)_\beta^\alpha (X)\bigg\downarrow\ \quad
& &
\ \quad\bigg\downarrow (\mu^{\overline F})_\beta^\alpha(X) \\
G^F_\beta(X) &
\stackrel{\displaystyle \nu_\beta(X)}{\hbox to 2cm{\rightarrowfill}}
& \ G_\beta^{\overline F}(X)
\end{array}
$$
is commutative for any manifold $X$.
  Then the natural transformation $\eta=\{\eta_\pi\}: F\to
\overline F$ given by the compositions
$$
\eta_\pi :
F(\pi)
\stackrel{\scriptstyle \kappa_\pi}{\hbox to 1cm{\rightarrowfill}}
T^{\mu^F}(\pi)
\stackrel{\scriptstyle \widetilde\nu_\pi}{\hbox to 1cm{\rightarrowfill}}
T^{\mu^{\overline F}} (\pi)
\stackrel{\scriptstyle \overline\kappa_\pi^{-1}}{\hbox to 1cm{\rightarrowfill}}
\overline F(\pi)
$$
is the unique natural transformation
$F\to \overline F$ such that
$\nu^\eta_\alpha=\nu_\alpha$, where $\nu^\eta_\alpha$
is defined by $(\ref{nu_eta})$.
}

\par
\medskip
\noindent
{\bf Definition.}
  We say that two bundle functors $F$ and
$\overline F$ are {\it equivalent} if there exists a natural
equivalence  $\eta: F\to \overline F$.
  We say that two inductive systems of natural transformations $\mu$
and
$\overline \mu$  are {\it equivalent} if there exists a family
$\nu=(\nu_\alpha)$ of natural transformations such that the diagram
(\ref{G-diag}) is commutative for any manifold $X$.
\par
\medskip
  The following proposition completes the proof of Theorem 3.1.
  It is proved just the same as Corollary 2.3  in~\cite{Mik}.
\par
\medskip
\noindent
{\bf Proposition 3.4.}
{\it
  The correspondence $F\to \mu^F$ induces
a bijection between the equivalence classes of product preserving
bundle functors on the category ${\cal FM}_\Lambda$ and the
equivalence classes of inductive systems of natural transformations
of product preserving bundle functors on the category ${\cal M}f$.
  The inverse bijection is induced by the correspondence
$\mu\to T^\mu$.}
\par

\bigskip

\begin{flushleft}
\sc
Geometry Department\\
Branch of Mathematics\\
Chebotarev Research Institute of Mathematics and Mechanics\\
Kazan State University\\
Universitetskaya, 17, Kazan, 420008\\
Russia\\[5pt]
{\it E-mail:} {\tt vadimjr@ksu.ru}
\end{flushleft}

\end{document}